\documentclass[12pt,reqno]{amsart}
\usepackage{amssymb}
\usepackage{a4wide}
\usepackage{amsfonts}    
\usepackage{amssymb,amsmath,amsthm,amscd}
\numberwithin{equation}{section}
\newtheorem{theorem}{Theorem}[section]
\newtheorem{proposition}[theorem]{Proposition}

\newtheorem{lemma}[theorem]{Lemma}

\theoremstyle{definition}

\theoremstyle{remark}
\newtheorem{remark}[theorem]{Remark}

\begin{document}

\title[ concentration of solutions]{concentration of solutions for a fourth order elliptic equation in $\mathbb{R}^N$}

\author{Zhongyuan Liu}
\address{Academy of Mathematics and Systems Science, Chinese Academy of Sciences, Beijing 100190,  P. R. China}
\email{liuzy@amss.ac.cn}


\begin{abstract}
 In this paper, we study the following
fourth order elliptic problem
$$
\Delta^2 u=(1+\varepsilon K(x)) u^{2^*-1},  \quad x\in  \mathbb{R}^N  \\
$$
where $2^*=\frac{2N}{N-4}$,$N\geq5$, $ \varepsilon>0$. We prove that
the existence of two peaks solutions for the above problem,  if
$K(x)$ has two critical points satisfying certain conditions,
provided $\varepsilon$ is small enough.
  \\[12pt]
 \emph{AMS 2000 subject classifications: Primary $35\mathrm{J}60$;\  Secondary   $35\mathrm{J}65$, $58\mathrm{E05}$ \newline
     }
\end{abstract}

\maketitle

\section{ Introduction and main results}

 In this paper, we
study the following nonlinear fourth order elliptic equation
\begin{equation}\label{11}
\left\{
\begin{array}{ll}
\Delta^2 u=(1+\varepsilon K(x)) u^{2^*-1}, u>0,  \quad x\in
\mathbb{R}^N\\
u\in \mathcal {D}^{2,2}(\mathbb{R}^N)
\end{array}
\right.
\end{equation}
where $ N\geq 5, \varepsilon>0$, $K(x)\in C^1(\mathbb{R}^N)\cap
L^\infty(\mathbb{R}^N)$ and  $\mathcal {D}^{2,2}(\mathbb{R}^N)$ be
the completion of $C_0^\infty(\mathbb{R}^N)$ with the respect to the
norm
$$\|u\|^2=\langle u,u\rangle, ~\text{where}  ~~ \langle
u,v\rangle=\int_{\mathbb{R}^N}\Delta u\Delta v.$$
$2^*=\frac{2N}{N-4}$ is the critical exponent of the embedding
$\mathcal {D}^{2,2}(\mathbb{R}^N)\hookrightarrow
L^{2^*}(\mathbb{R}^N).$ A Fourier transformation argument yields
$\mathcal {D}^{2,2}(\mathbb{R}^N)=\mathcal
{D}_0^{2,2}(\mathbb{R}^N)$.

In the past few years, there has been many study on concentration of
solutions for second-order elliptic equations with critical Sobolev
exponent; See e.g.\cite{APY1,APY2,AP,BP,CNY,CY,H,MP,R1,R,Y} and the
references therein. Recently, some researches have been developed on
the existence of peak solutions of fourth order elliptic equations
involving critical exponent,  see for example
\cite{BE,BEH,BS,CG,G,PZ}.

Define the Euler-Lagrange functional $I$ corresponding to \eqref{11}
as follows
$$
I_\varepsilon(u)=\frac{1}{2}\int_{\mathbb{R}^N} |\Delta
u|^2-\frac{1}{2^*}\int_{\mathbb{R}^N}(1+\varepsilon K(x))|u|^{2^*}.
$$
 One of the main features of problem \eqref{11} is the lack of
compactness, i.e. the functional $I_\varepsilon$ does not satisfy
the Palais-Smale condition. Such a fact follows from the
noncompactness of the embedding of $\mathcal
{D}^{2,2}(\mathbb{R}^N)\hookrightarrow L^{2^*}(\mathbb{R}^N)$ and
the unboundness of the domain $\mathbb{R}^N$. In this article, we
use a construction method to obtain peak solutions for \eqref{11}.
Precisely, we extend the argument employed in \cite{CNY} to the
framework of such higher order equations. To do so, we first take
advantage of a type of Lyapunov-Schmidit reduction to transform the
problem of finding critical points for the functional
$I_\varepsilon$ into one of finding critical
 points of a functional defined on finite dimensional domain. Then
 we construct a suitable bounded domain with finite dimension in
 which the associated variational problem, by topological degree argument, can have a critical point.
In our proof, to obtain a fine analysis on the energy of the
functional $I_\varepsilon$, we perform a careful expansion on
$I_\varepsilon$ by make full use of the precise  computation of the
contribution of function $K(x)$ to its critical points. Moreover, we
have to prove the positivity of the critical points obtained by our
process. It is well-known that such a proof, in general,  is quite
difficult for higher order equations.

To state the main result,  we need to introduce some notations and
assumptions.

Consider the equation
\begin{equation}\label{12}
\Delta^2 u=|u|^{2^*-2}u,~~u>0,~x\in\mathbb{R}^N, ~u\in \mathcal
{D}^{2,2}(\mathbb{R}^N).
\end{equation}
It has been proved in \cite{L} that the following function, for
$y\in\mathbb{R}^N$ and $ \lambda>0$,
$$
U_{y,\lambda}(x)=C_N\frac{\lambda^{\frac{N-4}{2}}}{(1+\lambda^2|x-y|^2)^{\frac{N-4}{2}}},~
C_N=[(N-4)(N-2)N(N+2)]^{\frac{N-4}{8}}
$$
solves \eqref{12} on $\mathbb{R}^N.$

Let $\Sigma$ denote the set consisting of all the critical points
$z$ of $K(x)$, satisfying (after a suitable rotation of the
coordinate system depending on $z$),
$$
K(x)=K(z)+\sum_{i=1}^{N}a_i|x_i-z_i|^\beta+O(|x-z|^{\beta+\sigma})
$$
for $x$ close to $z$, where $a_i, \beta$ and $\sigma$ are some
constants depending only on  $z$, $a_i\neq0$ for $i=1,2,\cdots,N$,
$\sum\limits_{i=1}^{N}a_i<0, \beta\in(1, N-4)$ and $\sigma\in(0,1)$.

We now state our main result of the paper:

\begin{theorem}\label{th11} \hspace{0.0cm}
Assume that $\Sigma$ contains at least two points. Then for each
$z^1, z^2\in\Sigma, z^1\neq z^2$, there exists an $\varepsilon_0>0$
such that \eqref{11} has a solution of the form
$u_\varepsilon=\sum\limits_{j=1}^{2}\alpha_{j,\varepsilon}U_{y_\varepsilon^j,\lambda_{j,\varepsilon}}
+v_\varepsilon$  if $\varepsilon\in(0, \varepsilon_0)$ with
$\alpha_{j,\varepsilon}\rightarrow1, y_\varepsilon^j\rightarrow z^j,
 \lambda_{j,\varepsilon}\rightarrow +\infty$ and $\int_{{\mathbb R}^N} |\Delta
 v_\varepsilon|^2 \rightarrow0 $ as $\varepsilon\rightarrow0$.
\end{theorem}

\begin{remark}
If the set $\Sigma$ contains $k (k\geq 2)$  points, it is east to
see that , by theorem1.1, there exist $\binom{k}{2}$ solutions for
problem\eqref{11}.
\end{remark}

\begin{remark}
Let $z^1,z^2,\cdots,z^k$ be $k$ local maximum points of $K(x)$
satisfying
$$
K(x)=K(z^i)+\sum_{j=1}^Na_{ij}|x_j-z_j^i|^\beta+O(|x-z^i|^{\beta+\sigma}),\quad
x\in B_\delta(z^i),
$$
where $\beta\in(1,N-4),a_{ij}\neq0$ and $\sum_{j=1}^Na_{ij}<0$ and
$\sigma>0$. Using the technique in the proof of Theorem 1.3 of
\cite{Y}, we can construct a $k$-peaked solution for \eqref{11},
such that there is exactly one local maximum point near each $z_i, i
= 1,2,\cdots, k$.
\end{remark}

This paper is organized as follows:  we first introduce some
notations and perform a finite-dimensional reduction in section 2,
and then use topological degree argument to prove theorem 1.1 in
section 3. In order that we can give a clear line of our framework,
we list all the proofs of the needed estimates in the appendix..

\section{Notations and the finite-dimensional reduction}

For
$y=(y^1,y^2,\cdots,y^k)\in\mathbb{R}^N\underbrace{\times\mathbb{R}^N\times\cdots\times}_k
\mathbb{R}^N$,
$\lambda=(\lambda_1,\lambda_2,\cdots,\lambda_k)\in\mathbb{R}^k$,
define
\begin{gather}\label{13}
E_{y,\lambda}^k=\{v\in \mathcal {D}^{2,2}(\mathbb{R}^N)|\langle
U_{y^j,\lambda_j},v\rangle=\left\langle\frac{\partial
U_{y^j,\lambda_j}}{\partial\lambda_j},v\right\rangle=\left\langle\frac{\partial
U_{y^j,\lambda_j}}{\partial y_i^j},v\right\rangle=0 \notag \\
\left.\text{for} j=1,2,\cdots,k, i=1,2,\cdots,N\right\}
\end{gather}

We look for solutions $u_\varepsilon$ of \eqref{11} of the following
form
\begin{equation}\label{14}
u_\varepsilon=\alpha_{1,\varepsilon}U_{y_\varepsilon^1,\lambda_{1,\varepsilon}}+\alpha_{2,\varepsilon}U_{y_\varepsilon^2,\lambda_{2,\varepsilon}}+v_\varepsilon
\end{equation}
with
$\alpha_\varepsilon=(\alpha_{1,\varepsilon},\alpha_{2,\varepsilon})\rightarrow
(1,1),
y_\varepsilon=(y_\varepsilon^1,y_\varepsilon^2)\rightarrow(z^1,z^2),
 \lambda_{j,\varepsilon}\rightarrow +\infty$ for $j=1,2$ and $v_\varepsilon\in E^2_{y_\varepsilon,\lambda_\varepsilon}$
  satisfying $\|v_\varepsilon\|\rightarrow 0 $ as $\varepsilon\rightarrow
  0$.

For each $z^1,z^2\in \Sigma, z^1\neq z^2$ and $\mu>0$, define
\[
\begin{split}
D_\mu=\left\{(y,\lambda)\Big|y=(y^1,y^2)\in
\overline{B_\mu(z^1)}\times\overline{B_\mu(z^2)},\lambda=(\lambda_1,\lambda_2)\in
\left(\frac{1}{\mu},+\infty\right)\times\left(\frac{1}{\mu},+\infty\right)\right\}.
\end{split}
\]

Define
\[
\begin{split}
M_\mu=\{(\alpha,y,\lambda,v)|\alpha=(\alpha_1,\alpha_2)\in\mathbb{R}^+\times\mathbb{R}^+,
(y,\lambda)\in D_\mu, v\in E_{y,\lambda}^2, \\
|\alpha_1-1|\leq\mu,|\alpha_2-1|\leq\mu,\|v\|\leq\mu
 \}
\end{split}
\]
and
$$
J_\varepsilon(\alpha,y,\lambda,v)=I_\varepsilon\left(\sum_{j=1}^2\alpha_jU_{y^j,\lambda_j}+v\right).
$$

It is well known now that for $\mu>0$ sufficiently small if
$(\alpha,y,\lambda, v)\in M_\mu$ is a critical point of
$J_\varepsilon$ in $M_\mu$, then
$u=\sum\limits_{j=1}^2\alpha_jU_{y^j,\lambda_j}+v$ is a critical
point of $I_\varepsilon$ in $\mathcal{D}^{2,2}({\mathbb{R}^N})$.The
fact that $(\alpha,y,\lambda, v)\in M_\mu$ is a critical point of
$J_\varepsilon$ in $M_\mu$ is equivalent to the fact that the
following equations are satisfied
\begin{align}\label{15}
\frac{\partial J_\varepsilon}{\partial
\alpha_j}&=0,~~j=1,2,\\\label{16} \frac{\partial
J_\varepsilon}{\partial
v}&=\sum\limits_{j=1}^2A_jU_{y^j,\lambda_j}+\sum\limits_{j=1}^2B_j\frac{\partial
U_{y^j,\lambda_j}}{\partial
\alpha_j}+\sum\limits_{j=2}^2\sum\limits_{i=1}^N
C_{ji}\frac{\partial U_{y^j,\lambda_j}}{\partial y_i^j} ,
\\\label{17} \frac{\partial J_\varepsilon}{\partial
y_i^j}&=B_j\left\langle\frac{\partial^2 U_{y^j,\lambda_j}}{\partial
\alpha_j\partial y_i^j}, v\right\rangle+\sum\limits_{l=1}^N
C_{jl}\left\langle\frac{\partial^2 U_{y^j,\lambda_j}}{\partial
y_l^j\partial y_i^j }, v\right\rangle, i=1,2,\cdots,N,
j=1,2,\\\label{18} \frac{\partial J_\varepsilon}{\partial
\lambda_j}&=B_j\left\langle\frac{\partial^2
U_{y^j,\lambda_j}}{\partial \alpha_j^2},
v\right\rangle+\sum\limits_{l=1}^N
C_{jl}\left\langle\frac{\partial^2 U_{y^j,\lambda_j}}{\partial
y_l^j\partial \lambda_j}, v\right\rangle,  j=1,2,
\end{align}
for some $A_j, B_j, C_{ji}\in\mathbb{R}, j=1,2, i=1,2,\cdots,N$.

As in \cite{CNY}, we first reduce the problem of finding a solution
for \eqref{11} to that of finding a critical point for a function
defined in a finite dimensional domain. Then we use topological
degree argument  to solve the latter problem. We next establish some
preliminary result. Throughout this paper we will let
$\varepsilon_{12}=\frac{1}{(\lambda_1\lambda_2)^{\frac{N-4}{2}}}$.

\begin{proposition}\label{p21} Suppose that $z^1, z^2\in\Sigma$ and
$(y,\lambda)\in D_\mu$. Then there exist $\varepsilon_0>0$ and
$\mu_0>0$ such that for $\varepsilon\in(0,\varepsilon_0)$ and
$\mu\in(0,\mu_0)$, there is a unique $C^1$ map $(y,\lambda)\in
D_\mu\rightarrow(\alpha_\varepsilon(y,\lambda),v_\varepsilon(y,\lambda))\in\mathbb{R}^2\times\mathcal{D}^{2,2}(\mathbb{R}^N)$
such that $v_\varepsilon\in E^2_{y,\lambda}$,
$(\alpha_\varepsilon,y,\lambda,v_\varepsilon)$ satisfies
\eqref{15},\eqref{18}. Furthermore,
$\alpha_\varepsilon=(\alpha_{1,\varepsilon},\alpha_{2,\varepsilon})$
and $v_\varepsilon$ satisfy the following estimate as
$\varepsilon\rightarrow0$
\begin{equation}\label{21}
\sum\limits_{j=1}^2|\alpha_{j,\varepsilon}-(1+\varepsilon
K(z^j))^{-\frac{N-4}{4}}|+\|v_\varepsilon
\|=O\left(\varepsilon\sum\limits_{j=1}^2\left(|y^j-z^j|^{\beta_j}+\frac{1}{\lambda_j^{\theta_j}}\right)+\varepsilon_{12}^{\frac{1}{2}+\tau}\right)
\end{equation}
where $\theta_j=\inf\{\beta_j,\frac{N+4}{2}\},\tau>0$ is some
constant.
\end{proposition}

\begin{proof}
Let $\hat{\alpha}=((1+\varepsilon K(z^1))^{-\frac{N-4}{8}},
(1+\varepsilon K(z^2))^{-\frac{N-4}{8}})$,
$\bar{\alpha}=\alpha-\hat{\alpha},
\omega=(\bar{\alpha},v)\in\mathbb{R}^2\times E^2_{y,\lambda}$.

As in \cite{CNY}(see also \cite{PZ}) we expand
$\hat{J}_\varepsilon(y,\lambda,\omega)=J_\varepsilon(\alpha,y,\lambda,v)$
at $\omega=0$ and obtain
\begin{equation}\label{22}
\hat{J}_\varepsilon(y,\lambda,\omega)=\hat{J}_\varepsilon(y,\lambda,0)+\langle
f_\varepsilon,\omega\rangle+\frac{1}{2}\langle
Q_\varepsilon\omega,\omega\rangle+R_\varepsilon(\omega),
\end{equation}
where $f_\varepsilon\in\mathbb{R}^2\times E^2_{y,\lambda}$ is a
linear form given by
\begin{equation}\label{23}
\begin{split}
\langle
f_\varepsilon,\omega\rangle=-\int_{\mathbb{R}^N}&(1+\varepsilon
K)\left(\sum\limits_{j=1}^2\hat{\alpha}_j
U_{y^j,\lambda_j}\right)^{2^*-1}v\\
&+\sum\limits_{k=1}^2\bar{\alpha}_k\Bigg(\left\langle\sum\limits_{j=1}^2\hat{\alpha}_jU_{y^j,\lambda_j},U_{y^k,\lambda_k}
\right\rangle \\
&-\int_{\mathbb{R}^N}(1+\varepsilon
K)\left(\sum\limits_{j=1}^2\hat{\alpha}_jU_{y^j,\lambda_j}\right)^{2^*-1}U_{y^k,\lambda_k}\Bigg),
\end{split}
\end{equation}

$Q_\varepsilon$ is a quadratic form on $\mathbb{R}^2\times
E^2_{y,\lambda}$ given by
\begin{equation}
\begin{split}
\langle Q_\varepsilon\omega,\omega\rangle &=\sum_ {k=1\atop
l=1}^2\bar{\alpha}_l\bar{\alpha}_k \Bigg(\left\langle
U_{y^l,\lambda_l,U_{y^k,\lambda_k}}\right\rangle\\&-(2^*-1)\int_{\mathbb{R}^N}(1+\varepsilon
K)\left(\sum_{j=1}^2\hat{\alpha}_jU_{y^j,\lambda_j}\right)^{2^*-2}U_{y^l,\lambda_l}U_{y^k,\lambda_k}\Bigg)\\
&+\|v\|^2-(2^*-1)\int_{\mathbb{R}^N}(1+\varepsilon
K)\left(\sum_{j=1}^2\hat{\alpha}_jU_{y^j,\lambda_j}\right)^{2^*-2}v^2\\
&-2(2^*-1)\sum_{k=1}^2\bar{\alpha}_k\int_{\mathbb{R}^N}(1+\varepsilon
K)\left(\sum_{j=1}^2\hat{\alpha}_jU_{y^j,\lambda_j}\right)^{2^*-2}U_{y^k,\lambda_k}v
\end{split}
\end{equation}
and $R_\varepsilon$ is the higher order term satiafying
\begin{equation}
D^{(i)}R_\varepsilon(\omega)=O(\|\omega\|^{2+\theta-i}),~~i=0,1,2
\end{equation}
where $\theta>0$ is some constant.

To show the existence of
$(\alpha_\varepsilon(y,\lambda),v_\varepsilon(y,\lambda))\in\mathbb{R}^2\times
E^2_{y,\lambda}$ so that
$(\alpha_\varepsilon,y,\lambda,v_\varepsilon)$ satisfy \eqref{16}
and \eqref{17}, it suffices to obtain
$\omega(y,\lambda)\in\mathbb{N}^2\times E^2_{y,\lambda}$ such that
$D\hat{J}_\varepsilon(y,\lambda,\omega)=0$ for each fixed
$(y,\lambda)\in D_\mu$, where $D$ stands for the derivative with
respect to $\omega$. $D\hat{J}_\varepsilon=0$ is equivalent to
\begin{equation}\label{26}
f_\varepsilon+Q_\varepsilon+DR_\varepsilon(\omega)=0.
\end{equation}

As in \cite{CNY}, it is not difficult by using Lemma
\ref{a4}-\ref{a6} to check that if $\mu>0,\varepsilon>0$ are small
enough, then for each $(y,\lambda)\in D_\mu, Q_\varepsilon$ is
invertible and there exists $C>0$, independent of $(y,\lambda)$,
such that $\|Q_\varepsilon^{-1}\|\leq C.$ So by the implicit
function theorem, we can prove that there is a $C^1-$map
$\omega(y,\lambda)$ satisfying \eqref{26}. Furthermore
\begin{equation}\label{27}
\|\omega\|\leq C\|f_\varepsilon\|.
\end{equation}

Applying Lemmas \ref{a1}-\ref{a3} in Appendix A to \eqref{23} we get
\begin{equation*}
|\langle
f_\varepsilon,\omega\rangle|=O\left(\varepsilon\sum\limits_{j=1}^2\left(|y^j-z^j|^{\beta_j}+
\frac{1}{\lambda_j^{\theta_j}}\right)+\varepsilon_{12}^{\frac{1}{2}+\tau}\right)\|\omega\|,
\end{equation*}
and consequently we obtain \eqref{21}. So we have completed the
proof of the proposition.
\end{proof}

It is worthwhile to point out that from Lagrange multiplier theorem,
there are $A_j,B_j,C_{ji} (j=1,2,i=1,2,\cdots,N)$ in $\mathbb{R}$
such that $(\alpha,y,\lambda,v)$ satisfy \eqref{15} and \eqref{16}.
So we only need to solve finite dimensional problem \eqref{17} and
\eqref{18}.

\section{Proof of main result}

In this section, we prove that for the $A_j,B_j,C_{ji}\in
\mathbb{R}$ obtained in the above section satisfying \eqref{15} and
\eqref{16} there exists $(\tilde{y},\tilde{\lambda})\in D_\mu$ such
that \eqref{17} and \eqref{18} are satisfied by
$(\tilde{\alpha},\tilde{y},\tilde{\lambda},\tilde{v})$. Firstly, we
give some estimates.
\begin{lemma}\label{l31}Let $(y,\lambda)\in D_\mu$, $(\alpha,v)$ be
obtained as in Proposition\ref{p21}. For $\mu>0$ and $\varepsilon>0$
small enough,we have for $k=1,2$
\[
\begin{split}
\frac{\partial
J_\varepsilon(\alpha,y,\lambda,v)}{\partial\lambda_k}= & -\frac{
C_{N,\beta_k}\varepsilon}{\lambda_k^{\beta_k+1}}\sum_{i=1}^Na_i^k+\frac{C_0\varepsilon_{12}}{\lambda_k|z^1-z^2|^{N-4}}
+O\left(\frac{\varepsilon\varepsilon_{12}}{\lambda_k}\right) \\
& +O\left(\frac{\varepsilon_{12}^{1+\tau_1}}{\lambda_k}\right)+
O\left(\frac{\varepsilon}{\lambda_k^{\beta_k}}|y^k-z^k|\right) \\
&+O\left(\frac{\varepsilon}{\lambda_k}\sum_{j=1}^2\left(\frac{1}{\lambda_j^{\beta_j+\sigma}}+|y^j-z^j|^{\beta_j+\sigma}\right)\right),
\end{split}
\]
where $C_{N,\beta_k}$ is a positive constant depending only on $N$£¬
and $\beta_k$, $C_0$ is a positive constant,
$\tau_1=\min\{\tau,\frac{1}{N-4}\}$.
\end{lemma}
\begin{proof}
Without loss of generality, we take $k=1$. Direct computations yield
\[
\begin{split}
&\frac{\partial
J_\varepsilon(\alpha,y,\lambda,v)}{\partial\lambda_1}=\alpha_1\Bigg[\left\langle\sum_{j=1}^2\alpha_jU_{y^j,\lambda_j}+v,
\frac{\partial U_{y^1,\lambda_1}}{\partial\lambda_1}\right\rangle\\
&\quad -\int_{\mathbb{R}^N}(1+\varepsilon
K)\left|\sum_{j=1}^2\alpha_jU_{y^j,\lambda_j}+v\right|^{2^*-2}\left(\sum_{j=1}^2\alpha_jU_{y^j,\lambda_j}+v\right)\frac{\partial
U_{y^1,\lambda_1}}{\partial\lambda_1}\Bigg].\\
\end{split}
\]
It is easy to check that
\[
\begin{split}
&\frac{\partial
J_\varepsilon(\alpha,y,\lambda,v)}{\partial\lambda_1}=\alpha_1\Bigg[\alpha_2\int_{\mathbb{R}^N}U_{y^2,\lambda_2}^{2^*-1}
\frac{\partial U_{y^1,\lambda_1}}{\partial\lambda_1}\\
&\quad -\int_{\mathbb{R}^N}(1+\varepsilon K)\frac{\partial
U_{y^1,\lambda_1}}{\partial\lambda_1}\Bigg(\sum_{j=1}^2\left(\alpha_j
U_{y^j,\lambda_j}\right)^{2^*-1}
+(2^*-1)\left(\alpha_1U_{y^1,\lambda_1}\right)^{2^*-2}(\alpha_2U_{y^2,\lambda_2})\Bigg)\Bigg]\\
&\quad -\alpha_1(2^*-1)\int_{\mathbb{R}^N}(1+\varepsilon
K)\frac{\partial
U_{y^1,\lambda_1}}{\partial\lambda_1}\left(\sum_{j=1}^2\alpha_jU_{y^j,\lambda_j}\right)^{2^*-2}v
+O\left(\frac{\varepsilon_{12}^{1+\tau_1}}{\lambda_1}\right)+\left(\frac{\|v\|^2}{\lambda}\right)\\
&\triangleq
I_1-I_2+O\left(\frac{\varepsilon_{12}^{1+\tau_1}}{\lambda_1}\right)+\left(\frac{\|v\|^2}{\lambda}\right).
\end{split}
\]
By Lemmas A.7, B.1-B.2, we can obtain
\[
\begin{split}
I_1&=\alpha_1(\alpha_2-\alpha_2^{2^*-1})\int_{\mathbb{R}^N}U_{y^2,\lambda_2}\frac{\partial
U_{y^1,\lambda_1}}{\partial\lambda_1}-(2^*-1)\alpha_1^{2^*-1}\alpha_2U_{y^1,\lambda_1}^{2^*-2}U_{y^2,\lambda_2}\\
&\quad
-\alpha_1^{2^*-1}\varepsilon\int_{\mathbb{R}^N}K\frac{\partial
U_{y^1,\lambda_1}}{\partial\lambda_1}U_{y^1,\lambda_1}^{2^*-1}+O\left(\frac{\varepsilon\varepsilon_{12}}{\lambda_1}\right)\\
&=-\alpha_1^{2^*-1}\varepsilon\int_{\mathbb{R}^N}K\frac{\partial
U_{y^1,\lambda_1}}{\partial\lambda_1}U_{y^1,\lambda_1}^{2^*-1}-(2^*-1)\alpha_1^{2^*-1}\alpha_2U_{y^1,\lambda_1}^{2^*-2}U_{y^2,\lambda_2}
+O\left(\frac{\varepsilon\varepsilon_{12}}{\lambda_1}\right)\\
&=-\frac{\varepsilon
C_{N,\beta_1}}{\lambda_1^{\beta+1}}\sum_{i=1}^Na_i^1+O\left(\frac{\varepsilon}{\lambda_1^{\beta_1}}|y^1-z^1|\right)
+O\left(\frac{\varepsilon}{\lambda_1^{\beta_1+1+\sigma}}\right)+O\left(\frac{\varepsilon}{\lambda_1}|y^1-z^1|^{\beta_1+\sigma}\right)\\
&\quad+O\left(\frac{C_0\varepsilon_{12}}{\lambda_1|z^1-z^2|^{N-4}}\right)+O\left(\frac{\varepsilon^2}{\lambda_1}\right)
+O\left(\frac{\varepsilon\varepsilon_{12}}{\lambda_1}\right)+O\left(\frac{\varepsilon_{12}^{\frac{N-2}{N-4}}}{\lambda_1}\right)
\end{split}
\]
and
\[
\begin{split}
I_2=O\left(\frac{\varepsilon_{12}^{\frac{1}{2}+\tau}}{\lambda_1}+\frac{\varepsilon}{\lambda_1}\sum_{j=1}^2
\left(\frac{1}{\lambda^{\theta_j}}+|y^j-z^j|^{\beta_j}\right)\right)\|v\|.
\end{split}
\]
Combining the above equalities, we can derive the conclusion.

\end{proof}

\begin{lemma}\label{l32} Under the same assumption as in
Lemma\ref{l31}, we have $(k,l=1,2,k\neq l)$
\[
\begin{split}
\frac{\partial J_\varepsilon(\alpha,y,\lambda,v)}{\partial y_i^k}= &
-D_{N,\beta_k}a_i^k\frac{
\varepsilon}{\lambda_k^{\beta_k-2}}(y_i^k-z_i^k)-C_1\frac{(y_i^k-y_i^l)}{(\lambda_1\lambda_2)^{\frac{N-4}{2}}}
+O\left(\lambda_k\varepsilon_{12}^{1+\tau_1}\right) \\
& +
O\left(\frac{\varepsilon}{\lambda_k^{\beta_k-1}}\lambda_k^2|y^k-z^k|^2\right)+O\left(\varepsilon\lambda_k\varepsilon_{12}\right)\\
&+O\left(\varepsilon\lambda_k\sum_{j=1}^2\left(\frac{1}{\lambda_j^{\beta_j+\sigma}}+|y^j-z^j|^{\beta_j+\sigma}\right)\right)
+O\left(\varepsilon_{12}^{\frac{N-1}{N-4}}\right),
\end{split}
\]
where $D_{N.\beta_k}$ is a positive constant depending only on $N$
and $\beta_k$, $C_1>0$ is a constant. $\tau_1$ and $\sigma$ are the
same as in Lemma\ref{l31}.
\end{lemma}

\begin{proof}
By using H\"{o}lder inequality and lemmas A.7, B.3-B.4 in the
Appendix, the calculation is similar  to the proof of Lemma
\ref{l31},  we omit the detail.
\end{proof}
\begin{lemma}\label{l33}
For $(y,\lambda)\in D_\mu$, let $(\alpha,v)\in\mathbb{R}^2\times
E_{y,\lambda}^2$ be obtained in Proposition \ref{p21}. Then
\[
\begin{split}
B_k&=O\left(\lambda_k^2\sum_{j=1}^2\left(\frac{\varepsilon}{\lambda_j^{\beta_j+1}}+\varepsilon|y^j-z^j|^{\beta_j+1}\right)\right)+O(\lambda_k\varepsilon_{12}),\\
C_{ki}&=O\left(\frac{\varepsilon}{\lambda_k^{\beta_k+1}}(\lambda_k|y^k-z^k|+\varepsilon_{12})\right)\\
&\quad+O\left(\frac{1}{\lambda_k}\sum_{j=1}^2\left(\frac{\varepsilon}{\lambda_j^{\beta_j-1+\hat{\sigma}}}+\varepsilon|y^j-z^j|^{\beta_j+\hat{\sigma}}\right)\right),
\end{split}
\]

where $\hat{\sigma}>0$ is some constant.
\end{lemma}
\begin{proof}
Let $\varphi\in\mathcal{D}^{2,2}_0(\mathbb{R}^N)$, then
\[
\begin{split}
\left\langle\frac{\partial J_\varepsilon}{\partial
v},\varphi\right\rangle=&\sum_{j=1}^2A_j\left\langle
U_{y^j,\lambda_j},\varphi\right\rangle+\sum_{j=1}^2B_j\left\langle\frac{\partial
U_{y^j,\lambda_j}}{\partial\lambda_j},\varphi\right\rangle\\
&+\sum_{j=1}^2\sum_{i=1}^NC_{ji}\left\langle\frac{\partial
U_{y^j,\lambda_j}}{\partial y_i^j},\varphi\right\rangle.
\end{split}
\]
Taking $\varphi=U_{y^k,\lambda_k}, \frac{\partial
U_{y^k,\lambda_k}}{\partial y_h^k}, \frac{\partial
U_{y^k,\lambda_k}}{\partial\lambda_k}, k=1,2,h=1,2,\cdots,N$
respectively and noting that
\[
\frac{\partial
J_\varepsilon}{\partial\lambda_k}=\left\langle\frac{\partial
J_\varepsilon}{\partial v},\frac{\partial
U_{y^k,\lambda_k}}{\partial\lambda_k}\right\rangle,\quad
\frac{\partial J_\varepsilon}{\partial
y_h^k}=\left\langle\frac{\partial J_\varepsilon}{\partial
v},\frac{\partial U_{y^k,\lambda_k}}{\partial y_h^k}\right\rangle,
\]
we obtain a quasi-diagonal linear system of equations of $A_j, B_j$
and $C_{ji}$, whose coefficients are given by
\[
\begin{split}
\bigg\langle U_{y^j,\lambda_j},U_{y^k,\lambda_k}\bigg\rangle=
\begin{cases}
E,\quad & j=k,\\
O(\varepsilon_{jk}),\quad & j\neq k,
\end{cases}
\end{split}
\]
\[
\begin{split}
\bigg\langle U_{y^j,\lambda_j},\frac{\partial
U_{y^k,\lambda_k}}{\lambda_k}\bigg\rangle=
\begin{cases}
0,\quad & j=k,\\
O(\frac{\varepsilon_{jk}}{\lambda_k}),\quad & j\neq k,
\end{cases}
\end{split}
\]
\[
\begin{split}
\bigg\langle U_{y^j,\lambda_j},\frac{\partial
U_{y^k,\lambda_k}}{y_h^k}\bigg\rangle=
\begin{cases}
0,\quad & j=k,\\
O(\lambda_k\varepsilon_{jk}),\quad & j\neq k,
\end{cases}
\end{split}
\]
\[
\begin{split}
\bigg\langle \frac{\partial
U_{y^j,\lambda_j}}{\lambda_j},\frac{\partial
U_{y^k,\lambda_k}}{\lambda_k}\bigg\rangle=
\begin{cases}
\frac{F}{\lambda_j^2},\quad & j=k,\\
O(\frac{\varepsilon_{jk}}{\lambda_j\lambda_k}),\quad & j\neq k,
\end{cases}
\end{split}
\]
\[
\begin{split}
\bigg\langle\frac{\partial U_{y^j,\lambda_j}}{\lambda_j},
\frac{\partial U_{y^k,\lambda_k}}{y_h^k}\bigg\rangle=
\begin{cases}
0,\quad & j=k,\\
O(\frac{\lambda_k\varepsilon_{jk}}{\lambda_j}), \quad & j\neq k,
\end{cases}
\end{split}
\]
\[
\begin{split}
\bigg\langle\frac{\partial U_{y^j,\lambda_j}}{y_l^j}, \frac{\partial
U_{y^k,\lambda_k}}{y_h^k}\bigg\rangle=
\begin{cases}
G\lambda_j^2\delta_{hl},\quad & j=k,\\
O(\lambda_j\lambda_k\varepsilon_{12}),\quad & j\neq k,
\end{cases}\\
\end{split}
\]
where $E,F,G$ are strictly positive constants depending only on $N$
and $\delta_{hl}$ is the Kronecker symbol.

Using the estimates of $\frac{\partial
J_\varepsilon}{\partial\lambda_k},\frac{\partial
J_\varepsilon}{\partial y_i^k}$ in Lemmas \ref{l31} and \ref{l32},
we can obtain the estimates of $B_k$ and $C_{ki}$.
\end{proof}

\begin{proof}[\bf{Proof of Theorem \ref{th11}}]
Let
$L_\varepsilon=\varepsilon^{{\beta_1\beta_2}/[\beta_1\beta_2-\frac{N-4}{2}(\beta_1+\beta_2)]}$.
To obtain the existence of solution $u_\varepsilon$ of the form
\eqref{14}, we only need to show that \eqref{17} and \eqref{18} are
satisfied by some $(y,\lambda)\in D_\mu$. We will show that for some
suitable $\delta>0,\gamma_1>0$ small and $\gamma_2>0$ large, there
exists $(y,\lambda)$ such that $(\lambda_1,\lambda_2)\in
[\gamma_1L_\varepsilon^{\beta_1^{-1}},\gamma_2L_\varepsilon^{\beta_1^{-1}}]\times[\gamma_1L_\varepsilon^{\beta_2^{-1}},\gamma_2L_\varepsilon^{\beta_2^{-1}}]$,
$y=(y^1,y^2)\in B_{\frac{\delta}{\lambda_1}}(z^1)\times
B_{\frac{\delta}{\lambda_2}}(z^2)$ together with the
$(\alpha(y,\lambda),y,\lambda,v(y,\lambda))$ satisfy \eqref{17} and
\eqref{18}.

From Proposition \ref{p21} and Lemmas \ref{l31}-\ref{l32}, we get
the following equivalent form of \eqref{17} and \eqref{18}:
\begin{equation}\label{31}
\begin{split}
\frac{\varepsilon}{\lambda_k^{\beta_k}}\lambda_k(y_i^k-z_i^k)&=
 O\left(\varepsilon\sum_{j=1}^2\left(\frac{1}{\lambda_j^{\beta_j+\sigma}}+|y^j-z^j|^{\beta_j+\sigma}\right)\right)\\
 &\quad+O\left(\frac{\varepsilon_{12}}{\lambda_k}\right), \quad\quad
 k=1,2,i=1,2,\cdots,N,
\end{split}
\end{equation}
\begin{equation}\label{32}
\begin{split}
\frac{\varepsilon}{\lambda_k^{\beta_k}}\sum_{i=1}^Na_i^k+\frac{d_k}{(\lambda_1\lambda_2)^{\frac{N-4}{2}}}=&O(\varepsilon\varepsilon_{12})
+O(\varepsilon_{12}^{1+2\tau_1})\\
&\quad+O\left(\varepsilon\sum_{j=1}^2\left(\frac{1}{\lambda_j^{\beta_j+\sigma}}+|y^j-z^j|^{\beta_j+\sigma}\right)\right)
\end{split}
\end{equation}
where $d_k,\tau_1$ are some positive constant.

Let
\[
\begin{split}
\lambda_1&=t_1L_\varepsilon^{\beta_1^{-1}},
\lambda_2=t_2L_\varepsilon^{\beta_2^{-1}}, t_1,t_2\in[\gamma_1,\gamma_2],\\
y^1-z^1&=\lambda_1^{-1}x^1,y^2-z^2=\lambda_2^{-1}x^2,x^1,x^2\in B_\delta(0).
\end{split}
\]
Then \eqref{31} and \eqref{32} can be rewritten in the following equivalent way
\begin{align}\label{33}
x^k&=o_\varepsilon(1) \quad k=1,2, \\\label{34}
t_k^{-\beta_k}+\left(\sum_{i=1}^Na_i^k\right)^{-1}d_k\left(t_1t_2\right)^{-\frac{N-4}{2}}&=o_\varepsilon(1)
\quad k=1,2.
\end{align}

Set
\[
\begin{split}
f(x^1,x^2)&=(x^1,x^2), \quad (x^1,x^2)\in\Omega_1\triangleq B_\delta(0)\times B_\delta(0)\\
g(t_1,t_2)&=(g_1(t_1,t_2),g_2(t_1,t_2)), \quad (t_1,t_2)\in\Omega_2\triangleq [\gamma_1,\gamma_2]\times[\gamma_1,\gamma_2]\\
g_k(t_1,t_2)&=\frac{1}{t_k^{\beta_k}}-\frac{m_k}{(t_1t_2)^{\frac{N-4}{2}}},\quad \text{where} ~~m_k=-\frac{d_k}{\sum_{i=1}^Na_i^k}>0.
\end{split}
\]
Then
\begin{equation}\label{35}
\deg(f,\Omega_1,0)=1.
\end{equation}

On the other hand, it is easy to see that $g=0$ has a unique solution $(t_1^*,t_2^*)$ in $[\gamma_1,\gamma_2]\times[\gamma_1,\gamma_2]$
if $\gamma_1$ is small and $\gamma_2>0$ is large. Furthermore,
\[
\begin{split}
\left.\frac{\partial g_1}{\partial t_1}\right|_{(t_1^*,t_2^*)}&=\frac{1}{t_1^*}\left(-\beta_1+\frac{N-4}{2}\right)
\frac{m_1}{(t_1^*t_2^*)^{\frac{N-4}{2}}},\\
\left.\frac{\partial g_1}{\partial t_2}\right|_{(t_1^*,t_2^*)}&=\frac{(N-4)m_1}{2t_2^*(t_1^*t_2^*)^{\frac{N-4}{2}}},\\
\end{split}
\]
and
\[
\begin{split}
\left.\frac{\partial g_2}{\partial t_1}\right|_{(t_1^*,t_2^*)}&=\frac{(N-4)m_2}{2t_1^*(t_1^*t_2^*)^{\frac{N-4}{2}}},\\
\left.\frac{\partial g_2}{\partial t_2}\right|_{(t_1^*,t_2^*)}&=\frac{1}{t_2^*}\left(-\beta_2+\frac{N-4}{2}\right)
\frac{m_2}{(t_1^*t_2^*)^{\frac{N-4}{2}}}.
\end{split}
\]
Then
\[
\text{Jac}\left.g\right|_{(t_1^*,t_2^*)}=\left(\beta_1\beta_2-\left(\beta_1+\beta_2\right)\frac{N-4}{2}\right)\frac{m_1m_2}{(t_1^*t_2^*)^{(N-3)}}<0.
\]
So $\deg(g,\Omega_2,0)=-1$. As a consequence
$$
\deg((f,g),\Omega_1\times\Omega_2,0)=\deg(f,\Omega_1,0)\times\deg(g,\Omega_2,0)=-1.
$$
Thus, \eqref{33}, \eqref{34} has a solution and  $
u_\varepsilon=\alpha_{1,\varepsilon}
U_{y_\varepsilon^1,\lambda_{1,\varepsilon}}+\alpha_{2,\varepsilon}U_{y_\varepsilon^2,\lambda_{2,\varepsilon}}+v_\varepsilon$
is a critical point of $I_\varepsilon$. By Proposition \ref{p31}, we
know $u_\varepsilon>0$  for  $\varepsilon> 0$ sufficiently small.
Thus the result follows.
\end{proof}

\begin{proposition}\label{p31}
Assume that $u_\varepsilon=\alpha_{1,\varepsilon}
U_{y_\varepsilon^1,\lambda_{1,\varepsilon}}+\alpha_{2,\varepsilon}U_{y_\varepsilon^2,\lambda_{2,\varepsilon}}+v_\varepsilon$
is a critical point of $I_\varepsilon$ and $v_\varepsilon$ satisfies \eqref{21}.Then for $\varepsilon > 0$ sufficiently small,
$u_\varepsilon > 0$.
\end{proposition}
\begin{proof}
we follow the idea in \cite{BEH},\cite{PZ} to prove the proposition. Set $u_\varepsilon^-=\max\{-u_\varepsilon,0\}$.\\
Let us introduce $\omega$ satisfying
\[
\Delta^2\omega=-(1+\varepsilon
K)\left(u_\varepsilon^-\right)^{\frac{N+4}{N-4}}
~~\text{in}~~\mathbb{R}^N,
\]
and note that $c_N|x|^{4-N}$ is the Green function of the operator
$\Delta^2$ in $\mathcal{D}_0^{2,2}(\mathbb{R}^N)$,\\
Thus, $\omega$ can be written as
$$
\omega=-c_N\int_{\mathbb{R}^N}|x-y|^{4-N}(1+\varepsilon
K)\left(u_\varepsilon^-\right)^{\frac{N+4}{N-4}}\text{d}x.
$$
It is easy to see that $\omega\leq0$  and
\[
\|\omega\|^2=-\int_{\mathbb{R}^N}(1+\varepsilon
K)\left(u_\varepsilon^-\right)^{\frac{N+4}{N-4}}\leq
C\|\omega\||u_\varepsilon^-|_{L^{2^*}}^{\frac{N+4}{N-4}}.
\]
Assume that $\|\omega\|\neq0$, then
$$
\|\omega\|\leq C|u_\varepsilon^-|_{L^{2^*}}^{\frac{N+4}{N-4}}.
$$

On the other hand
\[
\begin{split}
\|\omega\|^2&=-\int_{\mathbb{R}^N}(1+\varepsilon
K)\left(u_\varepsilon^-\right)^{\frac{N+4}{N-4}}\omega\\ &
\geq-\int_{\{u_\varepsilon\leq0\}}(1+\varepsilon
K)\left(u_\varepsilon^-\right)^{\frac{N+4}{N-4}}\omega+\int_{\{u_\varepsilon\geq0\}}(1+\varepsilon
K)\left(u_\varepsilon^+\right)^{\frac{N+4}{N-4}}\omega\\
&=\int_{\mathbb{R}^N}(1+\varepsilon
K)|u_\varepsilon|^{\frac{8}{N-4}}u_\varepsilon\omega=\int_{\mathbb{R}^N}\omega\Delta^2u_\varepsilon
=\int_{\mathbb{R}^N}u_\varepsilon\Delta^2\omega \\
&=-\int_{\mathbb{R}^N}(1+\varepsilon
K)\left(u_\varepsilon^-\right)^{\frac{N+4}{N-4}}u_\varepsilon=\int_{\mathbb{R}^N}(1+\varepsilon
K)(u_\varepsilon^-)^{2^*} \\
&\geq c^\prime|u_\varepsilon^-|^{2^*}_{L^{2^*}}.
\end{split}
\]
Thus,
$$
c^\prime|u_\varepsilon^-|^{2^*}_{L^{2^*}}\leq\|\omega\|^2\leq
c|u_\varepsilon^-|^{\frac{2(N+4)}{N-4}}_{L^{2^*}}.
$$
It is obvious that
$|u_\varepsilon^-|^{2^*}_{L^{2^*}}\leq|v_\varepsilon|_{L^{2^*}}$, So
for $\varepsilon$ sufficiently small, $u_\varepsilon^-\equiv0$,
which implies that $\omega\equiv0$ and we get a contradiction. As a
result, $\omega\equiv0$ and $u_\varepsilon^-\equiv0$. Therefore,
$u_\varepsilon>0$ and we complete the proof.
\end{proof}

\begin{appendix}
\section{}
In this appendix, we give some estimates and result that used in
Section 2.
\begin{lemma}\label{a1}\ For any $(y,\lambda)\in D_\mu$ and $v\in
E^2_{y,\lambda}$, there exists $\tau>0$ such that
\begin{equation*}
\int_{\mathbb{R}^N}K\left(\sum\limits_{j=1}^2\alpha_jU_{y^j,\lambda_j}\right)^{2^*-1}v=
O\left(\sum\limits_{j=1}^2\left(|y^j-z^j|^{\beta_j}+
\frac{1}{\lambda_j^{\theta_j}}\right)+\varepsilon_{12}^{\frac{1}{2}+\tau}\right)\|v\|.
\end{equation*}
 where $\tau>0$ and $\theta_j=\inf\{\beta_j,\frac{N+4}{2}\}.$
\end{lemma}
\begin{proof}
Use the inequality
\begin{gather}\label{01}
||a+b|^p-a^p-b^p|\leq
\begin{cases}
Ca^{p/2}b^{p/2},\quad 1<p<2,\\
Cab^{p-1}+Ca^{p-1}b,\quad p\geq2.
\end{cases}
\end{gather}
We have
\[
\begin{split}
\displaystyle\int_{\mathbb{R}^N}& K\left(\sum\limits_{j=1}^2\alpha_jU_{y^j,\lambda_j}\right)^{2^*-1}v\\
&=\displaystyle\int_{\mathbb{R}^N}K\sum_{j=1}^2\left(\alpha_jU_{y^j,\lambda_j}\right)^{2^*-1}v+
\left\{
\begin{array}{ll}
O\left(\sum\limits_{i\neq
j}\displaystyle\int_{\mathbb{R}^N}U_{y^i,\lambda_i}^{\frac{2^*-1}{2}}U_{y^j,\lambda_j}^{\frac{2^*-1}{2}}|v|\right),\quad
2<2^*<3, \\
O\left(\sum\limits_{i\neq
j}\displaystyle\int_{\mathbb{R}^N}U_{y^i,\lambda_i}^{2^*-2}U_{y^j,\lambda_j}|v|\right),\quad
2^*\geq3
\end{array}\right.\\
&=\sum\limits_{j=1}^2\alpha_j^{2^*-1}\int_{\mathbb{R}^N}K(x)U_{y^j,\lambda_j}^{2^*-1}v+
\begin{cases}
O\left(\sum\limits_{i\neq
j}\left(\displaystyle\int_{\mathbb{R}^N}U_{y^i,\lambda_i}^{\frac{2^*}{2}}U_{y^j,\lambda_j}^{\frac{2^*}{2}}\right)^{\frac{N+4}{2N}}\|v\|\right),\quad
2<2^*<3, \\
O\left(\sum\limits_{i\neq
j}\displaystyle\int_{\mathbb{R}^N}\left(U_{y^i,\lambda_i}^{\frac{16N}{(N+4)(N-4)}}U_{y^j,\lambda_j}^{\frac{2N}{N+4}}\right)^{\frac{N+4}{2N}}\|v\|\right),\quad
2^*\geq3
\end{cases}\\
&=\sum\limits_{j=1}^2\alpha_j^{2^*-1}\displaystyle\int_{\mathbb{R}^N}K(x)U_{y_j,\lambda_j}^{2^*-1}v+O\left(\varepsilon_{12}^{{1\over2}
+\tau}|v\|\right)\\
&=\sum\limits_{j=1}^2\alpha_j^{2^*-1}\displaystyle\int_{\mathbb{R}^N}(K(x)-K(z_j))U_{y_j,\lambda_j}^{2^*-1}v+O\left(\varepsilon_{12}^{{1\over2}
+\tau}|v\|\right)\\
&=O\left(\sum\limits_{j=1}^2\left(|y^j-z^j|^{\beta_j}+
\frac{1}{\lambda_j^{\theta_j}}\right)+\varepsilon_{12}^{\frac{1}{2}+\tau}\right)\|v\|.
\end{split}
\]
\end{proof}

\begin{lemma}\label{a2}\ There exists $\tau>0$ such that  for any $(y,\lambda)\in D_\mu$ and $v\in
E^2_{y,\lambda}$, we have
\begin{equation*}
\int_{\mathbb{R}^N}\left(\sum\limits_{j=1}^2\alpha_jU_{y^j,\lambda_j}\right)^{2^*-1}v=
O\left(\varepsilon_{12}^{\frac{1}{2}+\tau}\right)\|v\|.
\end{equation*}
\end{lemma}
\begin{proof}
By inequality\eqref{01}, we get
\[
\begin{split}
\int_{\mathbb{R}^N}&\left(\sum\limits_{j=1}^2\alpha_jU_{y^j,\lambda_j}\right)^{2^*-1}v\\
&=\sum_{j=1}^2\alpha_j^{2^*-1}\int_{\mathbb{R}^N}U_{y^j,\lambda_j}^{2^*-1}v+
\left\{
\begin{array}{ll}
O\left(\sum\limits_{i\neq
j}\displaystyle\int_{\mathbb{R}^N}U_{y^i,\lambda_i}^{2^*-1}U_{y^j,\lambda_j}^{\frac{2^*-1}{2}}|v|\right),\quad
2<2^*<3, \\
O\left(\sum\limits_{i\neq
j}\displaystyle\int_{\mathbb{R}^N}U_{y^i,\lambda_i}^{2^*-2}U_{y^j,\lambda_j}|v|\right),\quad
2^*\geq3
\end{array}\right.\\
&=O\left(\varepsilon_{12}^{\frac{1}{2}+\tau}\right)\|v\|
\end{split}
\]
since $\mu$ is small and $\lambda_1,\lambda_2\geq\frac{1}{\mu}$.
\end{proof}

\begin{lemma}\label{a3}\ Suppose $(y,\lambda)\in D_\mu$. Then for $\mu$ small enough we have
\begin{equation*}
\begin{split}
\left\langle\sum\limits_{j=1}^2\hat{\alpha}_jU_{y^j,\lambda_j},U_{y^k,\lambda_k}\right\rangle
&- \int_{\mathbb{R}^N}(1+\varepsilon
K)\left(\sum\limits_{j=1}^2\hat{\alpha}_jU_{y^j,\lambda_j}\right)^{2^*-1}U_{y^k,\lambda_k}\\
&=O\left(\varepsilon\left(\frac
{1}{\lambda_k^{\beta_k}}+|y^k-z^k|^{\beta_k}\right)+O(\varepsilon_{12})\right).
\end{split}
\end{equation*}
\end{lemma}

\begin{proof}
\[
\begin{split}
&\left\langle
\sum\limits_{j=1}^2\hat{\alpha}_jU_{y^j,\lambda_j},U_{y^k,\lambda_k}\right\rangle
-\int_{\mathbb{R}^N}(1+\varepsilon
K)\left(\sum\limits_{j=1}^2\hat{\alpha}_jU_{y^j,\lambda_j}\right)^{2^*-1}U_{y^k,\lambda_k}\\
&=\hat{\alpha}_k\left\langle
U_{y^k,\lambda_k},U_{y^k,\lambda_k}\right\rangle-\hat{\alpha}_k^{2^*-1}\int_{\mathbb{R}^N}(1+\varepsilon
K)U_{y^k,\lambda_k}^{2^*}+O(\varepsilon_{12})\\
&=\hat{\alpha}_k\left(\int_{\mathbb{R}^N}U_{y^k,\lambda_k}^{2^*}-\int_{\mathbb{R}^N}\frac{1+\varepsilon
K}{1+\varepsilon K(z^k)}U^{2*}_{y_k,\lambda_k}\right)+O(\varepsilon_{12})\\
&=\frac{\varepsilon}{(1+\varepsilon
K(z^k))^{\frac{N+4}{8}}}\int_{\mathbb{R}^N}\left(K(x)-K(z^k)\right)U^{2^*}_{y^k,\lambda_k}+O(\varepsilon_{12})\\
&=O(\varepsilon)\int_{\mathbb{R}^N}\left|K\left(\frac{z}{\lambda_k}+y^k\right)-K(z^k)\right|\frac{1}{(1+|z|^2)^N}+O(\varepsilon_{12})\\
&=O\left(\varepsilon\left(\frac
{1}{\lambda_k^{\beta_k}}+|y^k-z^k|^{\beta_k}\right)+O(\varepsilon_{12})\right).
\end{split}
\]
\end{proof}

\begin{lemma}\label{a4}\ Suppose $(y,\lambda)\in D_\mu$, $\mu$ and $\varepsilon$ small enough,
then there exists $\delta>0$ such that for all $v\in
E^2_{y,\lambda}$ we have
\begin{equation*}
\|v\|^2-(2^*-1)\int_{\mathbb{R}^N}(1+\varepsilon
K)\left(\sum\limits_{j=1}^2\hat{\alpha}_jU_{y^j,\lambda_j}\right)^{2^*-2}v^2\geq
\delta\|v\|^2.
\end{equation*}
\end{lemma}
\begin{proof}
It follows from Proposition 2 in \cite{PZ} that there exists a
constant $\delta_0>0$ such that
\[
\int_{\mathbb{R}^N}|\Delta
v|^2-(2^*-1)\int_{\mathbb{R}^N}\left(\sum_{j=1}^2\hat{\alpha}_jU_{y^j,\lambda_j}\right)^{2^*-2}v^2\geq\delta_0\int_{\mathbb{R}^N}|\Delta
v|^2
\]
for all $v\in E^2_{y,\lambda}$. Notice that
\[
\begin{split}
\Bigg|&\int_{\mathbb{R}^N}(1+\varepsilon
K)\left(\sum_{j=1}^2\hat{\alpha}_jU_{y^j,\lambda_j}\right)^{2^*-2}v^2-
\int_{\mathbb{R}^N}\left(\sum_{j=1}^2\hat{\alpha}_jU_{y^j,\lambda_j}\right)^{2^*-2}v^2\Bigg|\\
&\leq
O(\varepsilon)\|v\|^2+O\left(\int_{\mathbb{R}^N}U_{y^1,\lambda_1}^{\frac{2^*-2}{2}}U_{y^2,\lambda_2}^{\frac{2^*-2}{2}}v^2\right)\\
&=O\left(\varepsilon+(\lambda_1\lambda_2)^{-\hat{\sigma}}\right)\|v\|^2
\end{split}
\]
for some $\hat{\sigma}>0$ and Lemma A.4 follows.
\end{proof}

\begin{lemma}\label{a5}\ Suppose $(y,\lambda)\in D_\mu$.
Then  for $\mu$ and $\varepsilon$ small enough we have
\[
\begin{split}
\langle
U_{y^l,\lambda_l},U_{y^k,\lambda_k}\rangle-(2^*-1)\int_{\mathbb{R}^N}(1+\varepsilon
K)\left(\sum_{j=1}^2\alpha_jU_{y^j,\lambda_j}\right)^{2^*-2}U_{y^l,\lambda_l}U_{y^k,\lambda_k}\\
=\begin{cases}
\big(1-(2^*-1)\alpha_l^{2^*-2}\big)A+O(\varepsilon_{12}^r)+O(\varepsilon),~
&\text{if}~~ k=l=1,2,\\
O(\varepsilon_{12}^r)+O(\varepsilon),~&\text{if}~~ k\neq l,k,l=1,2,
\end{cases}
\end{split}
\]
for some $r>0,A>0.$
\end{lemma}
\begin{proof}
\[
\begin{split}
&\langle
U_{y^l,\lambda_l},U_{y^k,\lambda_k}\rangle-(2^*-1)\int_{\mathbb{R}^N}(1+\varepsilon
K)\left(\sum_{j=1}^2\alpha_jU_{y^j,\lambda_j}\right)^{2^*-2}U_{y^l,\lambda_l}U_{y^k,\lambda_k}\\
&=\int_{\mathbb{R}^N}U_{y^l,\lambda_l}^{2^*-1}U_{y^k,\lambda_k}-(2^*-1)\alpha_l^{2^*-2}\int_{\mathbb{R}^N}U_{y^l,\lambda_l}^{2^*-1}U_{y^k,\lambda_k}\\
&~~~~\quad+O\left(\int_{\mathbb{R}^N}U_{y^1,\lambda_1}^{2^*-2}U_{y^2,\lambda_2}^2+U_{y^1,\lambda_1}^{2^*-1}U_{y^2,\lambda_2}
U_{y^2,\lambda_2}^{2^*-2}U_{y^1,\lambda_1}^2+U_{y^2,\lambda_2}^{2^*-1}U_{y^1,\lambda_1}\right)+O(\varepsilon)\\
&=[1-(2^*-1)\alpha_l^{2^*-2}]\int_{\mathbb{R}^N}U_{y^l,\lambda_l}^{2^*-1}U_{y^k,\lambda_k}+O(\varepsilon_{12}^r)+O(\varepsilon).
\end{split}
\]
Let
$A=\int_{\mathbb{R}^N}U_{y^k,\lambda_k}^{2^*}=\int_{\mathbb{R}^N}U_{0,1}^{2^*}$,
we get the conclusion of Lemma A.5.
\end{proof}

\begin{lemma}\label{a6}\ Suppose $(y,\lambda)\in D_\mu, v\in E^2_{y,\lambda}$.
If $\mu$ and $\varepsilon$ small, then
\[
\begin{split}
&\int_{\mathbb{R}^N}(1+\varepsilon
K)\left(\sum_{j=1}^2\alpha_jU_{y^j,\lambda_j}\right)^{2^*-2}U_{y^k,\lambda_k}v\\
&=O\left(\varepsilon_{12}^{\frac1
2+\tau}+\varepsilon\sum_{j=1}^2\left(\frac{1}{\lambda_j^{\theta_j}}+|y^j-z^j|^{\beta_j}\right)\right)\|v\|
\end{split}
\]
 where $\tau>0$ and $\theta_j=\inf\{\beta_j,\frac{N+4}{2}\}.$
\end{lemma}
\begin{proof}
By using H\"{o}lder inequality and our assumption on $K(x)$, we have
\[
\begin{split}
\int_{\mathbb{R}^N}&(1+\varepsilon
K)\left(\sum_{j=1}^2\alpha_jU_{y^j,\lambda_j}\right)^{2^*-2}U_{y^k,\lambda_k}v\\
&=\alpha_k^{2^*-2}\int_{\mathbb{R}^N}U_{y^k,\lambda_k}^{2^*-1}v+\alpha_k^{2^*-2}\varepsilon\int_{\mathbb{R}^N}K(x)U_{y^k,\lambda_k}^{2^*-1}v\\
&\quad +\int_{\mathbb{R}^N}(1+\varepsilon
K)\left(\left(\sum_{j=1}^2\alpha_jU_{y^j,\lambda_j}\right)^{2^*-2}-\left(\alpha_kU_{y^k,\lambda_k}\right)^{2^*-2}\right)U_{y^k,\lambda_k}v\\
&=O\left(\varepsilon_{12}^{\frac1
2+\tau}+\varepsilon\sum_{j=1}^2\left(\frac{1}{\lambda_j^{\theta_j}}+|y^j-z^j|^{\beta_j}\right)\right)\|v\|.
\end{split}
\]
\end{proof}

\begin{lemma}\label{a7}\ Suppose $(y,\lambda)\in D_\mu, v\in E^2_{y,\lambda}$.
If $\mu$ and $\varepsilon$ small, then
\[
\begin{split}
\int_{\mathbb{R}^N}&(1+\varepsilon
K)\left(\sum_{j=1}^2\alpha_jU_{y^j,\lambda_j}\right)^{2^*-2}\frac{\partial
U_{y^k,\lambda_k}}{\partial
\lambda_k} v\\
&=O\left(\frac{\varepsilon_{12}^{\frac1
2+\tau}}{\lambda_k}+\frac{\varepsilon}{\lambda_k}\sum_{j=1}^2\left(\frac{1}{\lambda_j^{\theta_j}}+|y^j-z^j|^{\beta_j}\right)\right)\|v\|
\end{split}
\]
and
\[
\begin{split}
\int_{\mathbb{R}^N}&(1+\varepsilon
K)\left(\sum_{j=1}^2\alpha_jU_{y^j,\lambda_j}\right)^{2^*-2}\frac{\partial
U_{y^k,\lambda_k}}{\partial
 y_i^k} v\\
&=O\left(\lambda_k\varepsilon_{12}^{\frac1
2+\tau}+\lambda_k\varepsilon\sum_{j=1}^2\left(\frac{1}{\lambda_j^{\theta_j}}+|y^j-z^j|^{\beta_j}\right)\right)\|v\|
\end{split}
\]

 where $\tau>0$ and $\theta_j=\inf\{\beta_j,\frac{N+4}{2}\}.$
\end{lemma}
\begin{proof}
The proof of this lemma is similar to that of Lemma A.6.
\end{proof}

\section{}

  The estimates given in this appendix will play the key
role in proving our main result in Section 3. The computation here
is very similar the one performed in \cite{B}, so we give a sketch
proof.

\begin{lemma}\label{b1}\ Suppose $(y,\lambda)\in D_\mu$.
Then for $\mu$ small, $k=1,2$,  we have
\[
\begin{split}
\int_{\mathbb{R}^N}KU_{y^k,\lambda_k}^{2^*-1}\frac{\partial
U_{y^k,\lambda_k}}{\partial\lambda_k} &=\frac{C_{N,\beta_k}}{\lambda_k^{\beta_k+1}}\sum_{i=1}^Na_i^k+O\left(\frac{1}{\lambda_k^{\beta_k}}|y^k-z^k|\right)\\
 &+O\left(\frac{1}{\lambda_k^{\beta_k+1+\sigma}}\right)+O\left(\frac{1}{\lambda_k}|y^k-z^k|^{\beta_k+\sigma}\right)
\end{split}
\]
where $C_{N,\beta_k}$ is a positive constant depending only on $N$
and $\beta_k$.
\end{lemma}
\begin{proof}
\[
\begin{split}
\int_{\mathbb{R}^N}& KU_{y^k,\lambda_k}^{2^*-1}\frac{\partial
U_{y^k,\lambda_k}}{\partial\lambda_k}\\
&=\int_{\{|x-z^k|\leq
r_0\}}\left(\sum_{i=1}^Na_i^k|x_i-z_i^k|^{\beta_k}+O\left(|x-z^k|^{\beta_k+\sigma}\right)\right)U_{y^k,\lambda_k}^{2^*-1}\frac{\partial
U_{y^k,\lambda_k}}{\partial\lambda_k}
+O\left(\frac{1}{\lambda_k^{N+1}}\right)\\
&=\frac{C_N^\prime}{\lambda_k^{\beta_k+1}}\int_{\{|x-\lambda_k(z^k-y^k)|<\lambda_kr_0\}}\sum_{i=1}^ka_i^k|x_i-\lambda_k(y_i^k-z_i^k)|^{\beta_k}
\frac{1-|x|^2}{(1+|x|^2)^{N+1}}\\
&\quad
+O\left(\frac{1}{\lambda_k^{\beta_k+1+\sigma}}\right)+O\left(\frac{1}{\lambda_k}|y^k-z^k|^{\beta_k+\sigma}\right)\\
&=\frac{C_N^\prime}{\lambda_k^{\beta_k+1}}\int_{\mathbb{R}^N}\sum_{i=1}^Na_i^k|x_i|^{\beta_k}\frac{1-|x|^2}{(1+|x|^2)^{N+1}}
\\
&\quad
+O\left(\frac{1}{\lambda_k^{\beta_k}}|y^k-z^k|\right)+O\left(\frac{1}{\lambda_k^{\beta_k+1+\sigma}}+\frac{1}{\lambda_k}|y^k-z^k|^{\beta_k+\sigma}\right)\\
&=\frac{C_N^\prime}{\lambda_k^{\beta_k+1}}\sum_{i=1}^Na_i^k\int_{\mathbb{R}^N}|x|^{\beta_k}\frac{1-|x|^2}{(1+|x|^2)^{N+1}}
+O\left(\frac{1}{\lambda_k^{\beta_k}}|y^k-z^k|\right)\\
&\quad+O\left(\frac{1}{\lambda_k^{\beta_k+1+\sigma}}\right)+O\left(\frac{1}{\lambda_k}|y^k-z^k|^{\beta_k+\sigma}\right)\\
&\triangleq\frac{C_{N,\beta_k}}{\lambda_k^{\beta_k+1}}\sum_{i=1}^Na_i^k+O\left(\frac{1}{\lambda_k^{\beta_k}}|y^k-z^k|\right)
 +O\left(\frac{1}{\lambda_k^{\beta_k+1+\sigma}}+\frac{1}{\lambda_k}|y^k-z^k|^{\beta_k+\sigma}\right).
\end{split}
\]
\end{proof}

\begin{lemma}\label{b2}\ Suppose $(y,\lambda)\in D_\mu$,$\mu$ small,
 We have for $k\neq l,k,l=1,2,$
\[
\begin{split}
\int_{\mathbb{R}^N}U_{y^k,\lambda_k}^{2^*-2}\frac{\partial
U_{y^k,\lambda_k}}{\partial\lambda_k}U_{y^l,\lambda_l}=-\frac{N-4}{2}C_0\frac{\varepsilon_{12}}{\lambda_k
|z^1-z^2|^{N-4}}+O\left(\frac{\varepsilon_{12}^{\frac{N-2}{N-4}}}{\lambda_k}\right)
\end{split}
\]
where $C_0>0$.
\end{lemma}
\begin{proof}
\[
\begin{split}
\int_{\mathbb{R}^N}&U_{y^k,\lambda_k}^{2^*-2}\frac{\partial
U_{y^k,\lambda_k}}{\partial\lambda_k}U_{y^l,\lambda_l}=\frac{1}{2^*-1}\int_{\mathbb{R}^N}\frac{\partial
U_{y^k,\lambda_k}^{2^*-1}}{\partial\lambda_k}U_{y^l,\lambda_l}\\
&=\frac{1}{2^*-1}\int_{\mathbb{R}^N}\Delta\frac{\partial
U_{y^k,\lambda_k}}{\partial\lambda_k}\Delta U_{y^l,\lambda_l}
=\frac{1}{2^*-1}\int_{\mathbb{R}^N}\frac{\partial
U_{y^k,\lambda_k}}{\partial\lambda_k}U_{y^l,\lambda_l}^{2^*-1}\\
&=-\frac{N-4}{2}C_0\frac{\varepsilon_{12}}{\lambda_k
|z^1-z^2|^{N-4}}+O\left(\frac{\varepsilon_{12}^{\frac{N-2}{N-4}}}{\lambda_k}\right)
\end{split}
\]
(Follow the exact same line in proving estimate (F16) in \cite{B}).
\end{proof}

\begin{lemma}\label{b3}\ Suppose that $(y,\lambda)\in D_\mu$ and
 $\mu$ small.  Then
\[
\begin{split}
\int_{\mathbb{R}^N}KU_{y^k,\lambda_k}^{2^*-1}\frac{\partial
U_{y^k,\lambda_k}}{\partial y_i^k} &=D_{N,\beta_k} a_i^k\frac{1}{\lambda_k^{\beta_k-1}}\lambda_k(y^k_i-z^k_i)+
O\left(\frac{1}{\lambda_k^{\beta_k-1}}\lambda_k^2|y^k-z^k|^2\right)\\
 &+O\left(\frac{1}{\lambda_k^{\beta_k-1+\sigma}}\right)+O\left(\lambda_k|y^k-z^k|^{\beta_k+\sigma}\right)
\end{split}
\]
where $D_{N,\beta_k}$ is a positive constant depending only on $N$
and $\beta_k$.
\end{lemma}
\begin{proof}
\[
\begin{split}
&\int_{\mathbb{R}^N} KU_{y^k,\lambda_k}^{2^*-1}\frac{\partial
U_{y^k,\lambda_k}}{\partial
y_i^k}=(N-4)D_N\int_{\mathbb{R}^N}K(x)U_{y^k,\lambda_k}^{2^*}\frac{\lambda_k^2(x_i-y_i^k)}{1+\lambda_k^2|x-y^k|^2}\\
&=(N-4)D_N\int_{\{|x-z^k|\leq
r_0\}}\left(\sum_{h=1}^Na_h^k|x_h-z_h^k|^{\beta_k}+O(|x-z^k|^{\beta_k+\sigma})\right)\\
& \quad \times
U_{y^k,\lambda_k}^{2^*}\frac{\lambda_k^2(x_i-y_i^k)}{1+\lambda_k^2|x-y^k|^2}+O\left(\frac{1}{\lambda_k^{N-1}}\right)\\
&=\frac{D_N^\prime}{\lambda_k^{\beta_k-1}}\int_{\mathbb{R}^N}\sum_{h=1}^Na_h^k\left[|x|^{\beta_k}+\beta_k|x|^{\beta_k-1}x_h
\lambda_k(y_h^k-z_h^k)\right]\frac{U_{0,1}^{2^*}x_i}{1+|x|^2}\\
& \quad
+O\left(\frac{1}{\lambda_k^{\beta_k-1}}\lambda_k^2|y^k-z_k|^2\right)+O\left(\frac{1}{\lambda_k^{\beta_k-1+\sigma}}\right)
+O\left(\lambda_k|y^k-z^k|^{\beta_k+\sigma}\right)\\
\end{split}
\]
\[
\begin{split}
&=\frac{D_N^\prime\beta_k
a_i^k}{\lambda_k^{\beta_k-1}}\int_{\mathbb{R}^N}\frac{|x|^{\beta_k}}{(1+|x|^2)^{N+1}}\cdot\lambda_k(y_i^k-z_i^k)
+O\left(\frac{1}{\lambda_k^{\beta_k-1}}\lambda_k^2|y^k-z_k|^2\right)\\
&\quad+O\left(\frac{1}{\lambda_k^{\beta_k-1+\sigma}}\right)
+O\left(\lambda_k|y^k-z^k|^{\beta_k+\sigma}\right)\\
&\triangleq D_{N,\beta_k}
a_i^k\frac{1}{\lambda_k^{\beta_k-1}}\lambda_k(y^k_i-z^k_i)+
O\left(\frac{1}{\lambda_k^{\beta_k-1}}\lambda_k^2|y^k-z^k|^2\right)\\
&\quad+O\left(\frac{1}{\lambda_k^{\beta_k-1+\sigma}}\right)+O\left(\lambda_k|y^k-z^k|^{\beta_k+\sigma}\right).
\end{split}
\]
\end{proof}

\begin{lemma}\label{b4}\ Suppose that $(y,\lambda)\in D_\mu$ and $\mu$ small,
Then we have for $k\neq l,k,l=1,2,$
\[
\begin{split}
\int_{\mathbb{R}^N}U_{y^k,\lambda_k}^{2^*-2}\frac{\partial
U_{y^k,\lambda_k}}{\partial
y_i^k}U_{y^l,\lambda_l}=C_1\lambda_1\lambda_2(y_i^k-y_i^l)\varepsilon_{12}^{\frac{N-2}{N-4}}+O\left(\varepsilon_{12}^{\frac{N-1}{N-4}}\right)
\end{split}
\]
where $C_1>0$ is a constant depending only on $N$.
\end{lemma}
\begin{proof}
Follow the exact same procedure in proving estimate (F20) in
\cite{B}.
\end{proof}

\end{appendix}


\end{document}